\documentclass [12pt, a4paper] {article}

\usepackage {amsmath, amssymb, amsthm, amsfonts, natbib}

\title {Categoricity in Quasiminimal Pregeometry Classes}
\author {Levon Haykazyan\footnote{This work was completed through support of the
University of Oxford Clarendon Fund Scholarship.}}

\DeclareMathOperator {\dom} {dom}
\DeclareMathOperator {\img} {img}
\DeclareMathOperator {\cl} {cl}
\DeclareMathOperator {\id} {id}
\DeclareMathOperator {\cf} {cf}
\newcommand {\parto} {\rightharpoonup}

\theoremstyle {definition}
\newtheorem {definition} {Definition}
\newtheorem {example} [definition] {Example}

\newtheorem {theorem} [definition] {Theorem}
\newtheorem {proposition} [definition] {Proposition}
\newtheorem {corollary} [definition] {Corollary}
\newtheorem* {claim} {Claim}
\theoremstyle {remark}

\begin {document}
\maketitle

\begin {abstract}
Quasiminimal pregeometry classes were introduces by \cite {zilber} to isolate
the model theoretical core of several interesting examples. He proves that a
quasiminimal pregeometry class satisfying an additional axiom, called
excellence, is categorical in all uncountable cardinalities. Recently \cite
{bhhkk} showed that excellence follows from the rest of axioms. In this paper we
present a direct proof of the categoricity result without using excellence.
\end {abstract}

\section {Introduction}

Quasiminimal pregeometry class is a non elementary class of structures
satisfying certain axioms. The notion was introduced by \cite {zilber} to give
canonical axiomatisations of pseudo-exponential fields (in \cite {zilberps}) and
other related analytic structures. The original definition of \cite {zilber}
had an additional axiom called excellence, which played a central role in
establishing categoricity in cardinalities above $\aleph_1$. (And hence the
original terminology of a quasiminimal excellent class.) The notion has evolved
through works of \cite {baldwin} and \cite {kirby}. In practice checking that
the excellence holds has been the most technically difficult part in
applications of the categoricity theorem. Some of the original proofs of
excellence had gaps, which have only recently been fixed.

Later \cite {bhhkk} showed that excellence is redundant in that it follows from
the rest of axioms. In this paper we present a direct proof of categoricity
result that bypasses excellence altogether. The main new idea is to look at
(partial) embeddings that preserve all $\mathcal L_{\omega_1, \omega}$ formulas
possibly using infinitely many parameters. We call them $\sigma$-embeddings.
Constructing $\sigma$-embeddings by a transfinite recursion presents additional
challenges.  Given an increasing chain $\langle f_\beta : \beta < \alpha
\rangle$ of $\sigma$-embeddings, their union $f = \bigcup_{\beta < \alpha}
f_\beta$ need not be a $\sigma$-embedding. The problem is that $f$ needs to
preserve formulas over infinitely (countably) many parameters and if
$\cf(\alpha) = \omega$ we cannot guarantee that all these parameters occur at an
earlier stage. In case of quasiminimal pregeometry classes the axiom of
$\aleph_0$-homogeneity over countable closed models provides a way around this
problem in certain situations.

The rest of the paper is organised as follows. The next section fixes the
notation and gives basic definitions. Section \ref {properties} establishes
infinitary analogues of some well-known elementary properties in first order
model theory. Then we prove the categoricity theorem and finish with some
concluding remarks. 

The author is grateful to Martin Bays, Jonathan Kirby and Boris Zilber for
helpful discussions and comments.

\section {Background}

Let $\mathcal L$ be a finitary language, that is $\mathcal L$ has constant,
functional and relation symbols of finite arity. For an infinite cardinal
$\kappa$ the formulas of ${\mathcal L}_{\kappa, \omega}$ are inductively defined
as follows.
\begin {itemize}
\item Atomic $\mathcal L$-formulas are ${\mathcal L}_{\kappa, \omega}$-formulas.
\item If $\Phi$ is a set of ${\mathcal L}_{\kappa, \omega}$-formulas and $|\Phi|
	< \kappa$, then $\bigwedge_{\phi \in \Phi} \phi$ and $\bigvee_{\phi \in
	\Phi} \phi$ are ${\mathcal L}_{\kappa, \omega}$-formulas.
\item If $\phi$ is an ${\mathcal L}_{\kappa, \omega}$-formula and $v$ is a
	variable, then $\lnot \phi$, $\forall v\phi$ and $\exists v\phi$ are
	${\mathcal L}_{\kappa, \omega}$-formulas.
\end {itemize}
In this notation the ordinary first-order language coincides with ${\mathcal
L}_{\omega, \omega}$.  Note that we do not require ${\mathcal L}_{\kappa,
\omega}$-formulas to have finitely many free variables. However, every
subformula of an ${\mathcal L}_{\kappa, \omega}$-sentence has finitely many free
variables. If $M$ is an ${\mathcal L}$-structure, $\phi$ is an ${\mathcal
L}_{\kappa, \omega}$-formula  and $\theta$ is a variable assignment, then $M
\models_\theta \phi$ is defined as usual. In particular
$$M \models_\theta \bigwedge_{\phi \in \Phi} \phi \text { if and only if } M
\models_\theta \phi \text { for all } \phi \in \Phi$$
and
$$M \models_\theta \bigvee_{\phi \in \Phi} \phi \text { if and only if } M
\models_\theta \phi \text { for some } \phi \in \Phi.$$
Two structures $M$ and $N$ are called equivalent in ${\mathcal L}_{\kappa,
\omega}$ (in symbols $M \equiv_{\kappa, \omega} N$) if they satisfy the same
${\mathcal L}_{\kappa, \omega}$-sentences.

A formula of ${\mathcal L}_{\infty, \omega}$ is a formula of ${\mathcal
L}_{\kappa, \omega}$ for some $\kappa$. The notation $M \equiv_{\infty, \omega}
N$ means that $M$ and $N$ satisfy the same ${\mathcal L}_{\infty,
\omega}$-sentences.

Let $M, N$ be $\mathcal L$-structures and $f$ be a (partial) function from $M$
to $N$. (We use the notation $f : M \parto N$ for partial functions). Then $f$
is called a {\em (partial) embedding} if it preserves quantifier-free formulas.
Note that in particular $f$ preserves the formula $x = x$ and hence $f$ is
injective. A bijective embedding is an isomorphism between $M$ and $N$. The
function $f$ is called a {\em (partial) elementary embedding} if it preserves
first-order formulas and a {\em (partial) $\sigma$-embedding} if it preserves
all ${\mathcal L}_{\omega_1, \omega}$-formulas (in general using infinitely many
parameters from the $\dom(f)$). In other words $f$ is a $\sigma$-embedding if
$\langle M, m \rangle_{m \in \dom(f)} \equiv_{\omega_1, \omega} \langle N, f(m)
\rangle_{m \in \dom(f)}$.

A {\em back-and-forth system} between $\mathcal L$-structures $M$ and $N$ is a
nonempty collection $F$ of partial embeddings such that the following two
conditions are satisfied:
\begin {itemize}
\item for all $f \in F$ and $a \in M$ there is $g \in F$ such that $f \subseteq
	g$ and $a \in \dom(g)$;
\item for all $f \in F$ and $b \in N$ there is $g \in F$ such that $f \subseteq
	g$ and $b \in \img(g)$.
\end {itemize}

The following characterisation of ${\mathcal L}_{\infty, \omega}$-equivalence is
due to \cite {karp}.
%for a modern proof see \cite {marker-inf}.

\begin {theorem}
Let $M$ and $N$ be ${\mathcal L}$-structures. Then $M \equiv_{\infty, \omega} N$
if and only if there is a back and forth system between them.
\end {theorem}

Now we introduce the notion of a quasiminimal pregeometry class from \cite
{zilber}. Our definition follows closely that of \cite {kirby} . For the
definition of pregeometry the reader can consult \cite {marker}.

\begin {definition}
Let $\mathcal L$ be a language. A {\em quasiminimal pregeometry} class is a
class $\mathcal C$ of pairs $\langle H, \cl_H \rangle$ where $H$ is an $\mathcal
L$-structure and $\cl_H : {\mathcal P}(H) \to {\mathcal P}(H)$ is a function
satisfying the following conditions.
\begin {itemize}
\item {\em Closure under isomorphisms}\newline
	If $\langle H, \cl_H\rangle \in \mathcal C$, $H'$ is an $\mathcal
	L$-structure and $f : H \to H'$ is an isomorphism, then $\langle H',
	\cl_{H'} \rangle \in \mathcal C$, where $\cl_{H'}$ is defined as
	$\cl_{H'}(X') = f(\cl_H(f^{-1}(X')))$ for $X' \subseteq H'$.

\item {\em Quantifier free theory}\newline
	If $\langle H, \cl_H \rangle, \langle H', \cl_H' \rangle \in \mathcal
	C$, then $H$ and $H'$ satisfy the same quantifier free sentences. In
	other words the empty function is a partial embedding between any two
	structures.

\item {\em Pregeometry}
\begin {itemize}
\item For each $\langle H, \cl_H\rangle \in \mathcal C$ the function $\cl_H$ is
	a pregeometry on $H$ and the closure of any finite set is countable.
\item If $\langle H, \cl_H \rangle \in \mathcal C$ and $X \subseteq H$, then
	$\cl_H(X)$ together with the restriction of $\cl_H$ is in $\mathcal C$.
\item If $\langle H, \cl_H \rangle, \langle H', \cl_{H'} \rangle \in \mathcal
	C$, $X \subseteq H$, $y \in H$ and $f : H \parto H'$ is a partial
	embedding defined on $X \cup \{y\}$, then $y \in \cl_H(X)$ if and only
	if $f(y) \in \cl_{H'}(f(X))$.
\end {itemize}
\item {\em $\aleph_0$-homogeneity over countable closed models}\newline
	Let $\langle H, \cl_H \rangle, \langle H', \cl_{H'} \rangle \in \mathcal
	C$, subsets $G \subseteq H, G' \subseteq H'$ be countable closed or
	empty and $g : G \to G'$ be an isomorphism.
\begin {itemize}
\item If $x \in H, x' \in H'$ are independent from $G$ and $G'$ respectively,
	then $g \cup \{\langle x, x'\rangle\}$ is a partial embedding.
\item If $g \cup f : H \parto H'$ is a partial embedding, $X = \dom(f)$ is
	finite and $y \in \cl_H(X \cup G)$, then there is $y' \in H'$ such that
	$g \cup f \cup \{\langle y, y'\rangle\}$ is a partial embedding.
\end {itemize}
\end {itemize}
\end {definition}

To illustrate the definition we give some examples of quasiminimal pregeometry
classes.

\begin {example}
The class of models of a strongly minimal first order theory together with
algebraic closure is a quasiminimal pregeometry class, provided that the closure
of the empty set is infinite. The later condition is needed to ensure that the
closure of any subset of a given model is a model of the theory itself. This can
be readily checked by Tarski-Vaught test, using the strong minimality of the
theory.

Let the language $\mathcal L$ contain just one binary relation $E$. Consider the
class of $\mathcal L$-structures, where $E$ is an equivalence relation and each
equivalence class is countable. Define the closure of $X$ to be the set of
elements equivalent to some $x \in X$. This class is a quasiminimal pregeometry
class. It can be realised as the class of models of an $\mathcal L(Q)$-sentence.

Finally we mention mathematically interesting non elementary classes of
pseudo-exponential fields of \cite {zilberps} and group covers of \cite
{zilbergc}.
\end {example}

A class satisfying the above conditions and an additional condition
referred to as excellence, is called a {\em quasiminimal excellent} class in
\cite {zilber, baldwin, kirby}. It is shown in these works that any two
structures in a quasiminimal excellent class of the same uncountable cardinality
are isomorphic.  As mentioned above \cite {bhhkk} showed that excellence follows
from other conditions. The terminology of a quasiminimal pregeometry class comes
from \cite {bhhkk}, where the countable dimensional structure is called a
quasiminimal pregeometry structure. Thus combining the results of \cite {zilber}
and \cite {bhhkk} we get that two structures in a quasiminimal pregeometry class
of the same uncountable cardinality are isomorphic. In this paper we present a
direct proof of this result.

\section {Properties of Structures in Quasiminimal Pregeometry Classes}
\label {properties}

Fix a quasiminimal pregeometry class $\mathcal C$. The closure operator is often
understood, so we will simply refer to $\mathcal C$ as a class of structures
instead of a class of pairs. Given a structure $H \in \mathcal C$ and a
substructure $G \subseteq H$ also in $\mathcal C$, we denote $G \preceq H$ the
fact that $G$ is closed in $H$.

\begin {proposition}
Let $H, H' \in \mathcal C$, $H \preceq H'$ and $X \subseteq H$. Then $\cl_H(X) =
\cl_{H'}(X)$.
\end {proposition}

\begin {proof}
Consider the identity embedding from $H$ to $H'$. For $y \in H$ we have $y \in
\cl_H(X)$ if and only if $y \in \cl_{H'}(X)$. Hence $\cl_H(X) = \cl_{H'}(X) \cap
H$. But $\cl_{H'}(X) \subseteq \cl_{H'}(H) = H$. Hence $\cl_H(X) = \cl_{H'}(X)$.
\end {proof}

In view of this, we will drop the subscript from closure operator whenever no
confusion arises. Let us prove some direct consequences of
$\aleph_0$-homogeneity.

\begin {proposition}
\label {alephik}
Let $H, H' \in \mathcal C$, subsets $G \subseteq H, G' \subseteq H'$ be
countable closed or empty, $g : G \to G'$ be an isomorphism and $g \cup f : H
\parto H'$ be a partial embedding with $X = \dom(f)$, $X' = \img(f)$ finite.
\begin {itemize}
\item The mapping $g \cup f$ extends to an isomorphism $\hat g : \cl(X \cup G)
	\to \cl(X' \cup G')$.
\item If $y \in H \setminus \cl(X \cup G)$ and $y' \in H' \setminus \cl(X' \cup
	G')$, then $g \cup f \cup \{\langle y, y'\rangle\}$ is a partial
	embedding.
\end {itemize}
\end {proposition}

\begin {proof}
For the first assertion note that by countable closure property both $\cl(X \cup
G)$ and $\cl(X' \cup G')$ are countable. Let $\langle a_n : n < \omega \rangle$
and $\langle b_n : n < \omega \rangle$ enumerate $\cl(X \cup G)$ and $\cl(X'
\cup G')$ respectively. Construct an increasing family $\langle f_n : n <
\omega \rangle$ of finite mappings such that $g \cup f_n$ is a partial embedding
as follows. Let $f_0 = f$. For odd $n$ pick the least $m$ such that $a_m$ is not
in the domain of $f_{n-1}$. By $\aleph_0$-homogeneity there is $b \in \cl(X'
\cup G')$ such that $g \cup f_{n-1} \cup \{\langle a_m, b \rangle\}$ is a
partial embedding.  Put $f_n = f_{n-1} \cup \{\langle a_m, b \rangle\}$. For
even $n$ do the other way around.  Then $\hat g = \bigcup_{n < \omega} f_n$ is
an isomorphism between $\cl(X \cup G)$ and $\cl(X' \cup G')$.

For the second assertion extend $g \cup f$ to an isomorphism $\hat g : \cl(X
\cup G) \to \cl(X' \cup G')$. Now $\hat g \cup \{\langle y, y'\rangle\}$ is a
partial embedding by $\aleph_0$-homogeneity. Hence $g \cup f \cup \{\langle y,
y' \rangle\}$ is also a partial embedding.
\end {proof}

Next we introduce $\sigma$-types and prove $\sigma$-saturation of uncountable
structures in $\mathcal C$.

\begin {definition}
Let $H \in \mathcal C$, $A \subseteq H$ and $\bar v$ be a finite tuple of
variables.  A $\sigma$-type $p$ (in $H$) over $A$ in variables $\bar v$ is a set
of ${\mathcal L}_{\omega_1, \omega}$-formulas with parameters from $A$ and free
variables among $\bar v$ such that every countable subset is consistent with
$H$. That is for every countable $\Phi \subseteq p$ we have $H \models \exists
\bar v \bigwedge_{\phi \in \Phi} \phi(\bar v)$. 
\end {definition}

If the length of tuple $\bar v$ is $n$, then we call $p$ an $n$-type.  We can
think semantically of $\sigma$-types as a family of ${\mathcal L}_{\omega_1,
\omega}$-definable subsets such that each countable subfamily has a nonempty
intersection. A $\sigma$-type $p$ is {\em complete} if for every ${\mathcal
L}_{\omega_1, \omega}$-formula $\phi(\bar v)$ we have either $\phi \in p$ of
$\lnot \phi \in p$. This corresponds to a $\sigma$-complete ultrafilter on the
$\sigma$-algebra of ${\mathcal L}_{\omega_1, \omega}$-definable subsets. A
$\sigma$-type $p$ is called {\em isolated} if there is a consistent ${\mathcal
L}_{\omega_1, \omega}$-formula $\psi(\bar x)$ such that $H \models \forall \bar
x(\psi(\bar x) \to \phi(\bar x))$ for all $\phi \in p$. A $\sigma$-type $p$ is
{\em realised} in $H$ if $\bigcap_{\phi \in p} \phi(H^n) \neq \emptyset$.
Clearly each isolated $\sigma$-type is realised.  Since ${\mathcal L}_{\omega_1,
\omega}$-definable sets are closed under countable intersections, we have the
following

\begin {proposition}
If a $\sigma$-type contains a formula defining a countable set, then it is
isolated and hence realised.
\end {proposition}

Now let us study $\sigma$-types in quasiminimal pregeometry structures.

\begin {proposition}
Let $H \in \mathcal C$ and $X \subseteq H$ be countable. Let $a, b \in H
\setminus \cl(X)$. Then $a$ and $b$ realise the same $\sigma$-type over $X$.
\end {proposition}

\begin {proof}
Let $G = \cl(X)$ and $g_0 = \id_G \cup \{\langle a, b \rangle\}$. By
$\aleph_0$-homogeneity $g_0$ is an embedding. Now consider the collection $F$ of
finite extensions of $g_0$ that are embeddings. It is not empty as $g_0 \in F$.
We claim that $F$ is a back-and-forth system. Indeed if $g \in F$ and $y \in H$,
then $g = \id_G \cup f$, where $f$ has finite domain. If $y \in \cl(\dom(g))$,
then use $\aleph_0$-homogeneity to extend $g$ to $y$. Otherwise there is $y' \in
H \setminus \cl(\img(g))$. By Proposition \ref {alephik} the map $g \cup
\{\langle y, y'\rangle \}$ is an embedding. Similarly, we can extend $g$ to an
embedding with $y$ in the image.

Now expand the language by adding constant symbols $c_g$ for each $g \in G$ and
an additional constant $c$. Let $H_1$ be an expansion of $H$ by interpreting
$c_g$ by $g$ and $c$ by $a$. Similarly let $H_2$ be an expansion of $H$ by
interpreting $c_g$ by $g$ and $c$ by $b$. By the above there is a back-and-forth
system between $H_1$ and $H_2$. Hence $H_1 \equiv_{\infty, \omega} H_2$. In
particular $H_1 \equiv_{\omega_1, \omega} H_2$. It follows that every sentence
using parameters from $G$ that is true on $a$ is also true on $b$. Hence $a$ and
$b$ realise the same $\sigma$-type over $G$. Since $X \subseteq G$, elements $a$
and $b$ also realise the same $\sigma$-type over $X$.
\end {proof}

The next corollary establishes the analogy between quasiminimal pregeometry
structures and minimality in first-order context. It is also the motivation
behind the term quasiminimality.

\begin {corollary}
Let $H \in \mathcal C$ and $\phi(v)$ be an ${\mathcal L}_{\omega_1, \omega}$
formula (possibly using parameters). Then $\phi(H)$ is either countable or
cocountable.
\end {corollary}

\begin {proof}
Suppose otherwise. Let $\bar c$ be the parameters used in $\phi$. Then $\bar c$
is countable. Since both $\phi(H)$ and $\lnot \phi(H)$ are uncountable, there
are $a \in \phi(H) \setminus \cl(\bar c)$ and $b \in \lnot \phi(H) \setminus
\cl(\bar c)$.  This contradicts the fact that $a$ and $b$ realise the same
$\sigma$-type over $\bar c$.
\end {proof}

And now we establish the analogy between the closure in quasiminimal pregeometry
classes and algebraic closure.

\begin {corollary}
Let $H \in \mathcal C$ be uncountable, and $X \subseteq H$ be a countable
subset.  Then $y \in \cl(X)$ if and only if it satisfies an ${\mathcal
L}_{\omega_1, \omega}$-formula that has countably many solutions.
\end {corollary}

\begin {proof}
Assume that $y$ satisfies $\phi(v)$ that has countably many solutions. Since
$\cl(X)$ is countable, there is $y' \in \lnot \phi(H) \setminus \cl(X)$. Now $y$
and $y'$ do not satisfy the same type over $X$. Hence $y \in \cl(X)$.

Conversely assume that $y \in \cl(X)$. Pick $y' \in H \setminus \cl(X)$. Since
the closure is determined by the language, the map $\id_X \cup \{\langle
y, y' \rangle\}$ is not an embedding. Hence there is a quantifier-free formula
$\phi(v)$ over $X$ satisfied by $y$ but not $y'$. Now $\lnot \phi(H)$ cannot by
countable (as it implies that $y' \in \cl(X)$). Hence $\phi(H)$ is countable.
\end {proof}

Next we introduce the infinitary analogue of saturation and prove this
property for uncountable structures in quasiminimal pregeometry classes.

\begin {definition}
A structure $H$ is called $\sigma$-saturated if for every $X \subset H$ with
$|X| < |H|$ every $\sigma$-type over $X$ is realised in $H$.
\end {definition}

\begin {proposition}
Let $H \in \mathcal C$ be uncountable. Then $H$ is $\sigma$-saturated.
\end {proposition}

\begin {proof}
Let $X \subset H$ be a subset with $|X| < |H|$ and let $p$ be a $\sigma$-type
over $X$ in $n$ variables. We prove by induction on $n$ that $p$ is realised in
$H$.

Let $n = 1$. Put $G = \cl(X)$. Then $|G| = |X| + \aleph_0 < |H|$. So there is $y
\in H \setminus G$. If $y$ realises $p$, then we are done. So assume the
opposite. Then there is a formula $\phi(v) \in p$ such that $H \models \lnot
\phi(y)$. Now since $y \not \in \cl(X)$, we have $\lnot \phi(H)$ is uncountable.
Hence $\phi(H)$ is countable. But then $p$ is isolated and hence realised in
$H$.

Assume the hypothesis for $n$. Let $p$ be an $n+1$-type. As before let $G =
\cl(X)$. We claim that there is $x \in H$ such that $q_x = \{\phi(\bar v, x) :
\phi(\bar v, w) \in p\}$ is a $\sigma$-type. Assume the opposite. Then for every
$x$ there is a countable subset $p_x \subset p$ such that
$$H \models \lnot \exists \bar v \bigwedge_{\phi \in p_x} \phi(\bar v, x).$$
Pick $y \in H \setminus G$ and let $Y$ be the set of parameters used in formulas
of $p_y$. Then $Y \subseteq X$ is countable and $y \not \in \cl(Y)$. But since
every two elements outside $\cl(Y)$ realise the same type over $Y$, for every $z
\in H \setminus \cl(Y)$ we have
$$H \models \lnot \exists \bar v \bigwedge_{\phi \in p_y} \phi(\bar v, z).$$
Now let $W = p_y \cup \bigcup_{x \in \cl(Y)} p_x$. Then $W$ is countable and  we
have that
$$H \models \lnot \exists w \exists \bar v \bigwedge_{\phi \in W} \phi(\bar v,
w).$$
This contradicts the fact that $p$ is a $\sigma$-type. Thus for some $x \in H$
we have that $q_x = \{\phi(\bar v, x) : \phi(\bar v, w) \in p\}$ is a
$\sigma$-type. By induction hypothesis $q_x$ is realised in $H$ and hence $p$ is
also realised in $H$.
\end {proof}

\section {The Categoricity Theorem}

In this section we prove that any two structures of the same uncountable
cardinality in a quasiminimal pregeometry class are isomorphic. Let $H, H' \in
\mathcal C$ be of the same uncountable cardinality. Since we can construct a
back-and-forth system between $H$ and $H'$, we have that $H \equiv_{\omega_1,
\omega} H'$. In other words the empty embedding is a $\sigma$-embedding. In
analogy with first-order case we would like to extend a partial
$\sigma$-embedding to map $H$ onto $H'$. By $\sigma$-saturation we can extend
any $\sigma$-embedding to any one element (and recursively to any finite number
of elements). At limit stages however, we need to take unions. But the union of
$\sigma$-embeddings may not be a $\sigma$-embedding. However, the union of
$\sigma$-embeddings is always an embedding and in some cases this is sufficient
to get a $\sigma$-embedding.

\begin {proposition}
\label {sigmaem}
Let $H, H' \in \mathcal C$ be uncountable, subsets $G \subset H, G' \subset H'$
be countable closed and let $g : G \to G'$ be an isomorphism. Then $g$ is a
partial $\sigma$-embedding between $H$ and $H'$.
\end {proposition}

\begin {proof}
By $\aleph_0$-homogeneity and Proposition \ref {alephik} the set of embeddings
between $H$ and $H'$ that are finite extensions of $g$ is a back-and-forth
system. Hence if we add constant symbols for $G$ in $H$ and for $G'$ in $H'$ the
resulting structures will be ${\mathcal L}_{\omega_1, \omega}$-equivalent.
Therefore $g$ is a $\sigma$-embedding.
\end {proof}

We can use this result to extend a $\sigma$-embeddings to the closure of its
domain provided the later is countable.

\begin {proposition}
\label {extension}
Let $H, H' \in \mathcal C$ be uncountable and let $g : H \parto H'$ be a partial
$\sigma$-embedding with $X = \dom(g)$, $X' = \img(g)$ countable. Then $g$
extends to a $\sigma$-embedding $\hat g : H \parto H'$ with $\dom(\hat g) =
\cl(X)$ and $\img(\hat g) = \cl(X')$.
\end {proposition}

\begin {proof}
Let $\cl(X) = \{a_n : n < \omega\}$ and $\cl(X') = \{a'_n : n < \omega\}$.
Construct an increasing sequence $f_0 \subseteq f_1 \subseteq f_2 ...$ of
$\sigma$-embeddings as follows. Let $f_0 = g$. For even $n$, pick the least $m$
not in the domain of $f_n$. Let $p$ be the $\sigma$-type of $a_m$ over
$\dom(f_n)$. Consider the $\sigma$-type $p' =\{\phi(x, f_n(\bar b)) : \phi(x,
\bar b) \in p\}$. The set $p'$ is a $\sigma$-type as $f_n$ is a
$\sigma$-embedding. By $\sigma$-saturation of $H'$, the type $p'$ is realised by
some $c \in H'$. Let $f_{n+1} = f_n \cup \{\langle a_m, c \rangle \}$. For odd
$n$ go the other direction.

Now take $\hat g = \bigcup_{n < \omega} f_n$. Then $\hat g$ is an embedding
between countable closed set $\cl(X)$ and $\cl(X')$. By Proposition \ref
{sigmaem} the embedding $\hat g$ is a $\sigma$-embedding.
\end {proof}

In particular every countable embedding that extends to the closure of its
domain must be a $\sigma$-embedding.

\begin {corollary}
\label {sigmaemcor}
Let $H, H' \in \mathcal C$ be uncountable, $G \subset H$, $G' \subset H'$ be
countable closed or empty, $g : G \to G'$ an isomorphism, $a \in H \setminus G$
and $a' \in H' \setminus G'$. Then $g \cup \{\langle a, a' \rangle \}$ is a
$\sigma$-embedding.
\end {corollary}

Now we can prove the main result of this paper. The main difference between our
approach and the existing literature is the focus on $\sigma$-embeddings. The
existing proofs of categoricity start with an ordinary embedding (i.e. a
function that preserves quantifier free formulas) and extend it to an
isomorphism between two structures of the same cardinality. At certain stages of
the construction one needs to extend an embedding with domain of a special form
to its closure. The condition of excellence is precisely the statement that this
is possible. However, if we have a $\sigma$-embedding at hand, then we can
always extend it to the closure of its domain by Proposition \ref {extension}
(provided the domain is countable). This is where we bypass the need for
excellence.

\begin {theorem}
Let $H, H' \in \mathcal C$ be uncountable, let $B, B'$ be bases of $H, H'$
respectively and let $g : B \to B'$ be a bijection. Then $g$ extends to an
isomorphism $\hat g : H \to H'$.
\end {theorem}

\begin {proof}
Let $B = \{b_\alpha : \alpha < \kappa\}$, $B' = \{b'_\alpha : \alpha < \kappa
\}$ and $g(b_\alpha) = b'_\alpha$. For $n < \omega$, let $G_n = \cl(\{b_m : n
\le m < \omega\})$ and $G'_n = \cl(\{b'_m : n \le m < \omega\})$. By
$\aleph_0$-homogeneity there is an isomorphism $f_0 : G_0 \to G'_0$ extending
$g$. By Proposition \ref {sigmaem} the embedding $f_0$ is a $\sigma$-embedding.

For each finite subset $X \subset B$ we construct a number $n_X$ and a
surjective $\sigma$-embedding $f_X : \cl(G_{n_X}X) \to \cl(G'_{n_X}X')$ that
extends $g$ and satisfies the following condition: whenever $X \subseteq Y$, we
have $n_X \le n_Y$ and $f_X|_{\cl(G_{n_Y}X)} = f_Y|_{\cl(G_{n_Y}X)}$. 

Assume that we have constructed such embeddings. Define $\hat g : H \to H'$ as
follows. For every $x \in H$, we have $x \in \cl(X)$ for some finite $X \subset
B$. Define $\hat g(x) = f_X(x)$. By the assumption on the embeddings, the result
does not depend on the choice of $X$. Now $\hat g$ is surjective. Indeed for $x'
\in H'$ pick finite $X' \subset B'$ such that $x' \in \cl(X')$. Let $X =
g^{-1}(X')$. Then $x' \in \img(f_X)$. Since $f_X$ is an embedding $f_X^{-1}(x')
\in \cl(X)$. Hence $\hat g(f_X^{-1}(x')) = x'$. Also if $\bar x \in H$ is a
finite tuple, choose $X \subset B$ a finite set such that $\bar x \in \cl(X)$.
Then $\hat g(\bar x) = f_X(\bar x)$, preserves quantifier free formulas. Thus
$\hat g$ is an isomorphism.

We now proceed to the construction of $\sigma$-embeddings $f_X$ by  a
well-founded induction on partial order of subsets of $B$. Our construction is a
modification of that of \cite {baldwin, kirby}. Take $n_\emptyset = 0$ and
$f_\emptyset = f_0$. If $X = \{b_\alpha\}$ is a singleton do the following. If
$\alpha < \omega$, then take $n_X = \alpha + 1$, otherwise take $n_X = 0$. Then
$f_0 |_{G_{n_X}} : G_{n_X} \to G'_{n_X}$ is an isomorphism. By Corollary \ref
{sigmaemcor}, the map $f_0 \cup \{\langle b_\alpha, b'_\alpha \rangle\}$ is a
$\sigma$-embedding. So by Proposition \ref {sigmaem} it extends to an
isomorphism $f_X : \cl(G_{n_X}X) \to \cl(G'_{n_X}X')$.

Let $X = \{x_1, ..., x_l\}$, $l > 1$ and we have already constructed $n_Y$ and
$f_Y$ for every $Y \subsetneq X$. Let $n = \max \{n_Y : Y \subsetneq X\}$. Then
each $g_Y = f_Y|_{\cl(G_nY)}$ is a $\sigma$-embedding of $\cl(G_nY)$ onto
$\cl(G'_nY')$, where $Y' = g(Y)$. Now if $x \in \dom(g_{Y_1}) \cap
\dom(g_{Y_2})$ for $Y_1, Y_2 \subsetneq X$, then $x \in \dom(g_{Y_1 \cap Y_2})$.
Thus if we define $g_X = \bigcup_{Y \subsetneq X} g_Y$, then $g_X$ is a well
defined function.

Let $Y_i = X \setminus \{x_i\}$ and $h_k = \bigcup_{i = 1}^k g_{Y_i}$. So $g_X =
h_l$. Then $h_k$ is a mapping from $\bigcup_{i = 1}^k \cl(G_nY_i)$ onto
$\bigcup_{i = 1}^k \cl(G'_nY'_i)$. 

\begin {claim}
For each $k = 1, ..., l$ there is $m_k \ge n$ such that $h_k
|_{\bigcup_{i=1}^k \cl(G_{m_k}Y_i)}$ is a $\sigma$-embedding.
\end {claim}

\begin {proof}
By induction on $k$. For $k = 1$, we have that $h_k = g_{Y_1}$, hence we can
just take $m_k = n$.

Assume the hypothesis for $k-1$. Then there is $m = m_{k-1} \ge n$ such that
$h_{k-1}|_{\bigcup_{i = 1}^{k-1} \cl(G_mY_i)}$ is a $\sigma$-embedding.
Take $m_k = m + 1$. Let $C_{k-1} = \bigcup_{i = 1}^{k-1} \cl(G_{m+1}Y_i)$.
We should show that $h_k|_{C_{k-1} \cup \cl(G_{m+1}Y_k)}$ is a
$\sigma$-embedding. Let $\bar a \in C_{k-1}$ and $\bar b \in \cl(G_{m+1}Y_k)$ be
(possibly infinite) tuples. Let $\theta$ be an automorphism of $\cl(G_mX)$ that
fixes $\cl(G_{m+1}Y_k)$ and swaps $x_k$ with $b_m$

Let $e$ be a $\sigma$-embedding from $\cl(G_mX)$ onto
$\cl(G'_mX')$ extending the $\sigma$-embedding $h_{k-1}|_{\bigcup_{i=1}^{k-1}
\cl(G_mY_i)}$. Let $\tau = e \theta^{-1} e^{-1} g_{Y_k} \theta$. The
automorphism $\theta$ takes $C_{k-1} \cup \cl(G_{m+1}Y_k)$ to $\cl(G_mY_k)$ and
$g_{Y_k}$ is defined on it (as $m \ge n$). Hence $\tau$ is well defined. All
these maps preserve ${\mathcal L}_{\omega_1, \omega}$ formulas. Hence so does
$\tau$.

Let $Y_{ik} = Y_i \cap Y_k = X \setminus \{x_i, x_k\}$. Since $\bar a \in
C_{k-1} = \bigcup_{i = 1}^{k-1} \cl(G_{m+1}Y_i)$, we have that
$$\theta(\bar a) \in \bigcup_{i = 1}^{k-1} \cl(G_mY_{ik}).$$
We have that $g_{Y_k}$ agrees with $g_{Y_i}$ on $\cl(G_mY_{ik})$. Hence
$g_{Y_k}$ agrees with $h_{k-1}$ on $\bigcup_{i = 1}^{k-1} \cl(G_mY_{ik})$. Also
$e$ agrees with $h_{k-1}$ on $\bigcup_{i = 1}^{k-1} \cl(G_mY_i)$ and $h_k$
extends $h_{k-1}$. Thus
$$\tau(\bar a) = e\theta^{-1}e^{-1}g_{Y_k}\theta(\bar a) =
h_{k-1}\theta^{-1}h_{k-1}^{-1}h_{k-1}\theta(\bar a) = h_{k-1}(\bar a) = h_k(\bar
a).$$

The tuple $\bar b \in \cl(G_{m+1}Y_k)$. Hence it is fixed by $\theta$. The
embeddings $g_{Y_k}$ and $e$ preserve closure, hence $e^{-1}g_{Y_k}(b) \in
\cl(G_{m+1}Y_k)$ is fixed by $\theta^{-1}$. Hence
$$\tau(\bar b) = e\theta^{-1}e^{-1}g_{Y_k}\theta(\bar b) = ee^{-1}g_{Y_k}(\bar
b) = g_{Y_k}(\bar b) = h_k(\bar b).$$

Thus if $\phi$ is any ${\mathcal L}_{\omega_1, \omega}$-formula, then
$$H \models \phi(\bar a, \bar b) \iff H' \models \phi(\tau(\bar a), \tau(\bar
b)) \iff H' \models \phi(h_k(\bar a), h_k(\bar b)).$$
Thus $h_k$ is a $\sigma$-embedding on $\bigcup_{i = 1}^k \cl(G_{m+1}Y_i)$. 
\end {proof}

Thus there is some $r$ such that $g_X = h_l$ is a $\sigma$-embedding on
$\bigcup_{i = 1}^l \cl(G_rY_i)$. Put $n_X = r$. Then $h_l |_{\bigcup_{i = 1}^l
\cl(G_rY_i)}$ extends to a $\sigma$-embedding $f_X$ from the closure of the
domain $\cl(\bigcup_{i = 1}^l \cl(G_rY_i)) = \cl(G_rX)$ onto $\cl(G'_rX')$ by
Proposition \ref {extension}.
\end {proof}

\section {Concluding Remarks}

Ever since its introduction by Shelah, excellence has been the key notion for
categoricity in non-elementary classes. What our methods show is that in some
very natural mathematical examples one can use infinitary logic instead. The
fact that $\sigma$-embeddings and associated infinitary notions occur in natural
mathematical contexts is remarkable in itself. This opens up a possibility of a
broader use of infinitary logic both in elementary and non-elementary setting.

\bibliographystyle {plainnat}
\bibliography {bib}

\end {document}